\documentclass[12pt]{article}
\usepackage[spanish,activeacute]{babel}
\usepackage{amsmath,amsfonts,amssymb,amstext,amsthm,amscd,mathrsfs,amsbsy}
\usepackage[utf8]{inputenc}
\usepackage{xcolor}
\usepackage{hyperref}
\usepackage{graphics,graphicx}

\newtheorem{teo}{Teorema}[section]

\newtheorem{coro}[teo]{Corolario}

\theoremstyle{definition}
\newtheorem{defi}[teo]{Definici\'on}
\newtheorem{exam}[teo]{Ejemplo}
\newtheorem{rem}[teo]{Observaci\'on}

\newcommand{\V}{\mathbf V}

\newcommand{\K}{\mathbb K}
\newcommand{\R}{\mathbb R}
\newcommand{\N}{\mathbb N}
\newcommand{\Z}{\mathbb Z}
\newcommand{\Q}{\mathbb Q}
\newcommand{\C}{\mathbb C}

\newcommand{\G}{\operatorname{Grass}_d}

\newcommand{\Po}{\K[x_1,\ldots,x_n]}

\newcommand{\Si}{\mbox{Sing}}

\begin{document}

\title{Crónica de un contraejemplo}
\author{Daniel Duarte}
\date{}
\maketitle

\begin{abstract}
En los años sesenta, John Nash propuso un método para resolver singularidades de variedades algebraicas. Cinco décadas de resultados alentadores no lograron evitar el triste desenlace: el método no funciona en general. Este artículo cuenta la historia del ascenso y la caída de la explosión de Nash.
\end{abstract}


\section*{Introducci\'on}

\noindent\textit{In these interludes of, as it were, enforced rationality, I did succeed in doing some respectable mathematical research.}

\hspace{4cm}$-$John Nash, autobiograf\'ia para el Premio Nobel.
\medskip

En geometría algebraica estudiamos objetos geométricos formados por las soluciones de sistemas de ecuaciones polinomiales. Estos objetos son conocidos como variedades algebraicas. Un problema fundamental en esta área es la resolución de singularidades que, dicho de manera informal, consiste en \textit{suavizar} puntos de cruce o cuspidales que suelen aparecer en una variedad. El problema ha sido resuelto parcialmente y la búsqueda de una solución completa es un tema de intensa investigación en la actualidad.

En la década de los sesenta, John F. Nash Jr. propuso un método para abordar el problema de resolución, hoy en día conocido como la explosión de Nash. Se trata de un proceso que remplaza puntos singulares por la tangencia en puntos vecinos. Ahora bien, debido a los problemas mentales que lo abrumaban, Nash no pudo desarrollar su método; le tocaría al mundo explorarlo, ponerlo en acción. Entre 1975 y 2025, matemáticas y matemáticos alrededor del mundo le dieron vida al método propuesto por Nash, obteniendo muchos resultados prometedores. Hasta el día en que un grupo de matemáticos latinoamericanos encontraron un contraejemplo.

En este artículo contamos la historia de la explosión de Nash. Tras introducir rápidamente el problema de resolución y el método de Nash, haremos un recorrido cronológico por los resultados obtenidos en esos cincuenta años, que podría resumirse como sigue:
\begin{itemize}
\item 1975-1977: los años de los primeros resultados.
\item 1977-1990: la época de oro, cuando reconocidas personalidades de la teoría de singularidades estudiaron el problema.
\item 1990-2010: tiempos de escasez, en los que la explosión de Nash cayó en un relativo olvido.
\item 2010-2012: el episodio tórico, en el que la explosión de Nash de variedades tóricas súbitamente se puso de moda.
\item 2020-2024: los años que marcan el estudio de la explosión de Nash sobre campos de característica positiva.
\item 2024: el año del contraejemplo.
\item 2025: donde comentamos los trabajos más recientes.
\end{itemize}

Aclaramos que el énfasis se hará en la narración de las diferentes etapas en la evolución del problema, evitando o manteniendo al mínimo los detalles técnicos propios de cada resultado. 

Sirva este artículo como un homenaje a la comunidad matemática que durante cinco décadas exploró un método tan fascinante como misterioso para la resolución de singularidades. La respuesta final representa la culminación del trabajo de todos esos años.


\section{El problema de la resolución de singularidades}

Empecemos definiendo el objeto básico de estudio. Un tratamiento detallado de este objeto, nociones relacionadas y sus propiedades básicas, se puede consultar en cualquier texto introductorio a la geometría algebraica (ver, por ejemplo, \cite{Hu}).

\begin{defi}
Sean $\K$ un campo y $S\subset\Po$ un subconjunto. Llamamos \textit{variedad algebraica definida por S} al conjunto
$$\V(S)=\{p\in \K^n\mid f(p)=0\mbox{ para todo }f\in S\}.$$
\end{defi}

Conocemos ejemplos de variedades algebraicas desde nuestra educación matemática elemental: rectas, parábolas, circunferencias, elipses e hipérbolas,  son todos ejemplos de conjuntos de puntos que anulan a algún polinomio. Veamos ahora un par de ejemplos menos familiares. 

\begin{figure}[ht]\label{curva nodal}
\begin{center}
\includegraphics[scale=0.35]{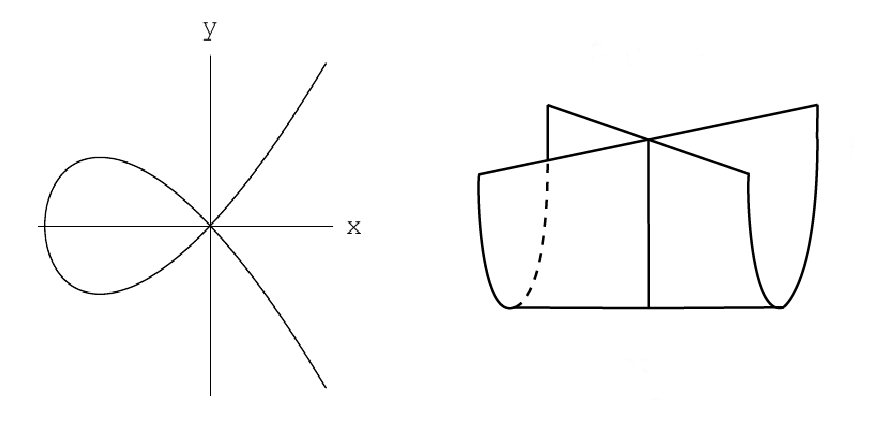}
\caption{Las variedades $\V(f)$ y $\V(g)$.}
\end{center}
\end{figure}

Sean $f(x,y)=x^3+x^2-y^2$ y $g(x,y,z)=x^2z-y^2$. Las variedades algebraicas definidas por $f$ y $g$ se ilustran en la figura \ref{curva nodal}.

Notemos una diferencia sustancial entre la curva $\V(f)$ y los dibujos que conocemos de la geometría analítica básica: el punto $(0,0)$ en la curva $\V(f)$ es un punto especial, no hay un punto que tenga un comportamiento parecido en una parábola o un círculo. Este ejemplo ilustra la noción de \textit{punto singular}. 

Sin entrar en definiciones formales, de manera intuitiva entendemos por puntos singulares de una variedad algebraica aquellos puntos en los que no está definido el espacio tangente. Denotamos al conjunto de puntos singulares de una variedad $X$ como $\Si(X)$. 

Ahora podemos enunciar el problema de resolución de singularidades: 
\medskip

\noindent\textit{Dada una variedad algebraica X, encontrar una variedad Y y una función polinomial sobreyectiva $\pi:Y\to X$ tal que:}
\begin{itemize}
\item[(i)] $\pi:Y\setminus\pi^{-1}(\Si(X))\to X\setminus\Si(X)$ \textit{es una biyección}.
\item[(ii)] $Y$ \textit{no tiene singularidades.}
\end{itemize}

El siguiente resultado, debido a Heisuke Hironaka y que lo hizo merecedor de una medalla Fields en 1970, da una respuesta afirmativa al problema sobre campos de característica cero (es decir, cualquier campo que contenga a $\Q$).

\begin{teo}\cite{H}
Toda variedad algebraica definida sobre un campo de característica cero admite una resolución de singularidades.
\end{teo}

En cambio, para campos de característica positiva (es decir, cualquier campo que contenga a $\Z/p\Z$ para algún número primo $p$), el problema sigue abierto y es un tema de intensa investigación en la actualidad.


\section{La explosi\'on de Nash}

En el Congreso Mundial de Psiquiatr\'ia celebrado en Madrid en 1996, John Nash declaró: 
\begin{center}
\textit{I had an idea which is referred to as Nash blowing up which I discussed \\ with an eminent mathematician named Hironaka}. 
\end{center}
La idea a la que se refiere y que ahora es com\'unmente conocida como la explosi\'on de Nash, tuvo lugar alrededor de 1963 durante una estancia que realiz\'o en el \textit{Institut for Advanced Study}, en Princeton (ver \cite[Chapter 42]{Na}).

Para dar un primer acercamiento a la explosión de Nash, consideremos nuevamente la curva $C=\V(x^3+x^2-y^2)\subset\R^2$ (ver figura \ref{curva nodal}). Para esta curva se cumple que $\Si(C)=\{(0,0)\}$. La idea consiste en construir una variedad $C^*$ usando los espacios tangentes en los puntos no singulares de $C$. Dado $p\in C$, el espacio tangente a $C$ en $p$ es denotado $T_pC$.

De entrada, para cada punto $p\in C\setminus\{(0,0)\}$, $T_pC$ es un subespacio vectorial de $\R^2$ de dimensión 1. Por el contrario, en el origen la curva no tiene un espacio tangente bien definido.

\begin{figure}[ht]\label{limites}
\begin{center}
\includegraphics[scale=0.5]{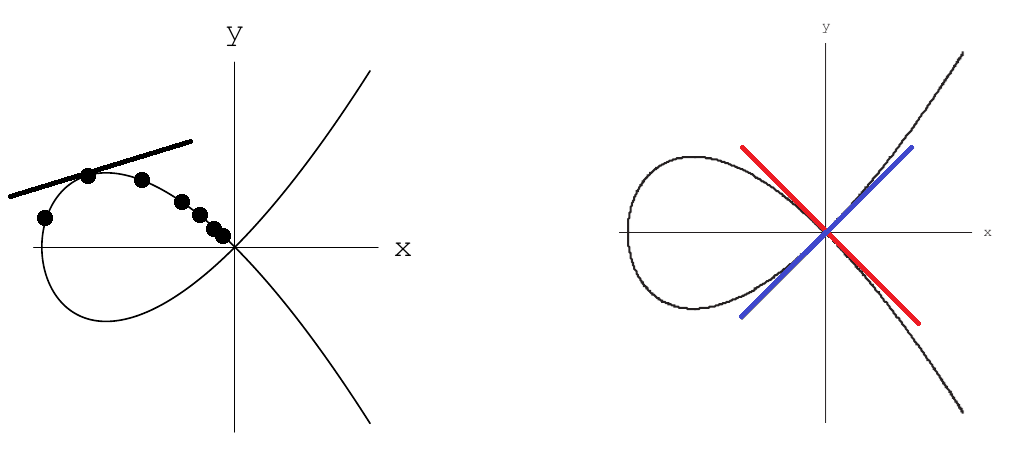}
\caption{Los l\'imites de espacios tangentes en $(0,0)$.}
\end{center}
\end{figure}

Y aquí viene la idea de Nash: asociar al punto singular todos los posibles l\'imites de espacios tangentes sobre sucesiones de puntos no singulares en $C$ que converjan a $(0,0)$. La figura \ref{limites} ilustra una sucesión de puntos que converge al punto singular y la recta tangente a uno de ellos. A un lado, los dos límites posibles de rectas tangentes, que corresponden a las líneas:
$$l_1=\{(x,y)\in\R^2|x-y=0\},\,\,\,\,\,l_2=\{(x,y)\in\R^2|x+y=0\}.$$


De esta manera, en el nuevo espacio $C^*$, el punto singular $p=(0,0)$ se convierte en dos parejas distintas: $(p,l_1)$ y $(p,l_2)$. ¡El punto singular se separó! (ver figura 3).
\begin{figure}[ht]\label{Nash de curva}
\begin{center}
\includegraphics[scale=0.6]{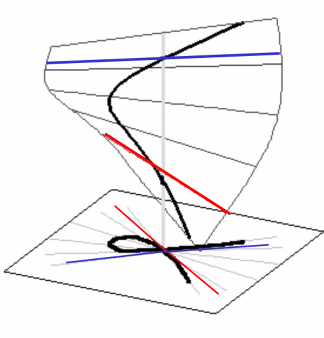}
\caption{La curva nodal desplegándose.}
\end{center}
\end{figure}

Enseguida procedemos a definir formalmente la explosi\'on de Nash de una variedad algebraica en $\K^n$. 

\begin{defi}
Sea $X\subset\K^n$ una variedad algebraica irreducible de dimensión $d$. Denotemos por $\G(\K^n)$ al conjunto de subespacios vectoriales de dimensi\'on $d$ en $\K^n$. Consideremos ahora la aplicaci\'on de Gauss:
\begin{align}
G:X\setminus\Si&(X)\to X\times\G(\K^n),\notag\\
&p\mapsto (p,T_pX).\notag
\end{align}
Sea $X^*$ la variedad algebraica más pequeña en $X\times\G(\K^n)$ que contiene a la imagen de $G$ (es decir, la cerradura de $Im(G)$ en la topología de Zariski). Sea $\nu:X^*\to X$ la restricci\'on  a $X^*$ de la proyecci\'on $X\times\G(\K^n)\to X$. La pareja $(X^*,\nu)$ es la \textit{explosión de Nash} de $X$. 
\end{defi}


Surge entonces la pregunta: si $X$ es una variedad singular, ?`$X^*$ lo sigue siendo? De ser el caso, repetimos el procedimiento. La pregunta que planteó Nash es la siguiente:
\begin{center}
\textit{?`Una iteraci\'on finita de la explosi\'on de Nash da lugar \\ a una variedad sin singularidades?}
\end{center}

De tener una respuesta afirmativa, este procedimiento daría lugar a un algoritmo canónico para resolver singularidades. 

\begin{rem}
Cabe mencionar que John Nash no publicó nada sobre el método que propuso. Al parecer esta pregunta se la hizo directamente a Hironaka en una carta privada (ver la introducción de \cite{Sp}). Por otro lado, aunque la pregunta se le atribuye usualmente a Nash, también aparece en el trabajo de J. G. Semple, diez años antes (ver \cite{S}). Algunos autores han incorporado a Semple denominando \textit{la explosión de Semple-Nash}.
\end{rem}

\section{La historia (1975-2025)}

Una vez establecidos los conceptos básicos y la pregunta a estudiar, podemos empezar a relatar la historia del problema. En las siguientes páginas haremos un recorrido por los resultados que se obtuvieron en cincuenta años sobre las propiedades de resolución de la explosión de Nash. 

\subsection{Los primeros resultados (1975-1977)}

En 1975 Augusto Nobile inaugura el estudio de la explosión de Nash con dos resultados que ya podemos considerar como clásicos: sobre campos de característica cero el método no es trivial, mientras que en característica positiva tristemente no funciona. 

\begin{teo}\cite[Theorem 1]{No}\label{Nob}
Sea $X$ una variedad algebraica sobre $\C$. La explosión de Nash $\nu:X^*\to X$ es un isomorfismo si y sólo si la variedad $X$ no tiene singularidades.
\end{teo}

Como corolario de este resultado se obtiene la primera respuesta afirmativa a la pregunta de Nash.

\begin{coro}\cite[Corollary 1]{No}\label{Nob curvas}
La explosión de Nash resuelve las singularidades de las curvas complejas.
\end{coro}

La pregunta no corrió la misma suerte sobre campos de característica prima, como Nobile ilustró con el siguiente ejemplo.

\begin{exam}\cite[Example 1]{No}\label{Nob ejem}
Sea $C=\V(x^3-y^2)\subset\K^2$ y asumir que la característica de $\K$ es 2. De entrada, Nobile muestra que calcular $C^*$ es equivalente a una operación llamada la explosión del ideal jacobiano de $C$, es decir, el ideal $\langle 3x^2,-2y\rangle$. 
Esto se denota usualmente como:
$$C^*\cong Bl_{\langle 3x^2,-2y \rangle}C.$$
Por otro lado, como la característica de $\K$ es 2, tenemos $\langle 3x^2,-2y \rangle=\langle 3x^2\rangle$. Una propiedad básica de las explosiones de ideales nos dice que explotar un ideal principal es un isomorfismo. Así:
$$C^*\cong Bl_{\langle 3x^2,-2y \rangle}C=Bl_{\langle 3x^2\rangle}C\cong C.$$

En conclusión: ¡$C^*\cong C$! La iteración del procedimiento no resolverá las singularidades de $C$.
\end{exam}

¡Ay de la explosión de Nash, que muere tan pronto como nace en característica positiva! Llena de pesar, la comunidad matemática la entierra y le pone flores. A partir de este momento toda la atención se centrará sobre el campo de los números complejos. 

Un resultado más se logra en estos años, obra de Vaho Rebassoo.

\begin{teo}\cite[Theorem 3.1]{R}
Sea $S=\V(x^py^q-z^r)\subset\C^3$, donde $p,q,r\in\N$ no tienen factor común. La iteración de la explosión de Nash resuelve las singularidades de $S$.
\end{teo}


\subsection{La época de oro (1977-1990)}

Los años que van de 1977 a 1990 serán la época de oro para la explosión de Nash. Tres distinguidas personalidades de la teoría de singularidades se darán cita para darle un espectacular impulso a este método.

Entre 1977 y 1982, Gerardo González-Sprinberg publica una serie de artículos en los que estudia una variante de la pregunta de Nash, que consiste en la iteración de la explosión de Nash seguida de la normalización\footnote{La normalización de una variedad algebraica es un proceso que elimina las singularidades en codimensión uno.}.
Nos referiremos a esta variante como la \textit{explosión de Nash normalizada} (ENN). González-Sprinberg introdujo técnicas que le permitieron estudiar la ENN con herramientas de geometría convexa, para ciertas familias de superficies. Este novedoso enfoque será clave para el descubrimiento del contraejemplo, casi cincuenta años después.

\begin{teo}\label{GS forever}
\begin{itemize} 
\item[(i)]\cite[\S 2, Théorème 1]{GS1} La ENN de una variedad tórica\footnote{Las variedades tóricas, grandes protagonistas de esta historia, son variedades algebraicas que tienen una descripción combinatoria en términos de semigrupos o conos poliedrales.} normal compleja coincide con la normalización de la explosión del ideal logarítmico jacobiano.
\item[(ii)]\cite[\S 2.3, Théorème]{GS1} La iteración de la ENN resuelve las singularidades de las superficies tóricas normales complejas.
\item[(iii)]\cite[\S 5]{GS2} La iteración de la ENN resuelve las singularidades de las superficies complejas ADE.
\end{itemize}
\end{teo}
¡La ENN inicia muy bien!

\begin{rem}\label{log jac 1}
\textit{Ideal logarítmico jacobiano (primera parte).} En este momento de la historia aparece el protagonista de la explosión de Nash en el contexto tórico: el ideal logarítmico jacobiano. El teorema \ref{GS forever} $(i)$ representa el fundamento teórico para describir de manera combinatoria la explosión de Nash normalizada de variedades tóricas complejas. Como iremos viendo en el relato, este ideal evolucionará hacia la versión más general posible, esto es, eliminando la condición de normalidad y para campos de característica arbitraria. Llegar hasta esa generalidad tomará más de cuatro décadas.
\end{rem}

En estos mismos años hace su aparición el mismísimo Heisuke Hironaka. Hacia 1979 Hironaka se interesa en la ENN y da nueva evidencia de que funciona. 

\begin{teo}\cite[Theorem]{Hi}
La iteración de la ENN da lugar a uniformización local a lo largo de valuaciones divisoriales\footnote{La uniformización local es una versión débil del problema de resolución, el objetivo es resolver solo un punto singular de la variedad que está determinado por una valuación.}, para cualquier variedad algebraica compleja.
\end{teo}

Ahora bien, Hironaka también demuestra que su teorema tiene la siguiente consecuencia en el caso de dimensión dos: la iteración de la ENN produce deliciosas superficies sandwich\footnote{Una superficie sandwich es una superficie que domina y es dominada por superficies no singulares.}. 
Convencido de que el método funciona en dimensión dos, Hironaka le pasa la estafeta a su estudiante Mark Spivakovsky para que concluya el caso de superficies. Y Mark Spivakovsky lo logra en 1990, en un esfuerzo descomunal que dio lugar a un artículo de ochenta páginas publicado en \textit{Annals of Mathematics}.

\begin{teo}\cite[Theorem 2.1]{Sp}
La iteración de la ENN resuelve las singularidades de las superficies complejas.
\end{teo}
La época de oro fue un periodo que inició con la pregunta de González-Sprinberg y que concluyó con el trabajo de Spivakovsky. Al respecto, en la introducción de \cite{Sp}, Spivakovsky nos obsequió una línea que unifica el esfuerzo de esos años:
\begin{center}
\textit{Also, our proof of desingularization of minimal singularities is a \\ generalization of, and inspired by, González-Sprinberg's\\ proof of the cyclic quotient case.}
\end{center}

\subsection{Tiempos de escasez (1990-2010)}

A años de abundancia le siguieron otros de escasez. Tras los espectaculares avances de las tres luminarias singularistas, vino una larga temporada de silencio. En los siguientes veinte años habrá pocas aportaciones a la pregunta de Nash. 

En 2003, un breve pero importante destello aparece en el apéndice de un trabajo de Monique Lejeune-Jalabert y Ana Reguera. Ahí se da una nueva demostración del teorema \ref{GS forever} $(i)$ (ver \cite[Proposition A.1]{LJ-R}). Esta nueva demostración representó el primer paso hacia la descripción combinatoria de la explosión de Nash de variedades tóricas sin la hipótesis de normalidad.

En 2006 tiene lugar una nueva demostración del caso de curvas planas, en manos de Monique Lejeune-Jalabert. A diferencia del corolario \ref{Nob curvas}, que se deduce del teorema \ref{Nob} a través de un argumento de noetherianidad, el enfoque de Lejeune-Jalabert muestra explícitamente el proceso de resolución.

\begin{teo}\cite[Theorem 3.1]{LJ}
La multiplicidad de un punto en la explosión de Nash de una curva plana coincide con la multiplicidad del punto en una transformación cuadrática de la curva. En particular, la iteración de la explosión de Nash resuelve las singularidades de la curva.
\end{teo}

Finalmente, en 2009 Wolfgang Ebelling y Sabir Gusein-Zade retoman la pregunta original y dan una nueva respuesta positiva.

\begin{teo}\cite[Section 1]{EGZ}
La explosión de Nash de una variedad determinantal genérica compleja es no singular.
\end{teo}

Aunque en su momento estos resultados no tuvieron mucho eco en la comunidad, en algunos años reaparecerían con toda su fuerza.

\subsection{El episodio tórico (2010-2012)}

Y de repente, ¡bum! Entre 2010 y 2012 la explosión de Nash de variedades tóricas se pone de moda. En este periodo cuatro grupos de matemáticos retomaron aquellas preguntas: Pedro González y Bernard Teissier, Dima Grigoriev y Pierre Milman, Atanas Atanasov, Christopher Lopez, Alexander Perry, Nicholas Proudfoot y Michael Thaddeus y, finalmente, el que escribe estas líneas, nos entregamos en cuerpo y alma al estudio de ambas preguntas en el caso de variedades tóricas. Aún con todo el esfuerzo y las diferentes perspectivas que cada equipo aportó, el caso tórico se resistió; nadie logró demostrarlo en su cabal generalidad. 

\begin{teo}
Los siguientes enunciados se refieren a variedades tóricas complejas.
\begin{itemize}
\item \cite[Section 6]{At} Implementación de la descripción combinatoria de la ENN en computadora. Se calculan miles de ejemplos en dimensión 3 y 4, en todos ellos se obtiene una resolución.
\item\cite[Theorem 2.1]{GM} Nueva demostración del teorema \ref{GS forever} $(ii)$. Además se dan cotas para la cantidad de iteraciones necesarias para llegar a la resolución.
\item \cite[Proposition 60]{GT} La explosión de Nash de una variedad tórica coincide con la explosión del ideal logarítmico jacobiano.
\item \cite[Theorem 109]{GT} Uniformización local a lo largo de valuaciones de rango maximal, en dimensión arbitraria.
\item \cite[Theorem 6.6]{D1} Uniformización local a lo largo de un conjunto denso de valuaciones, en dimensión dos.
\end{itemize}
\end{teo}

\begin{rem}\label{log jac 2}
\textit{Ideal logarítmico jacobiano (segunda parte).} En sus orígenes, la definición de variedad tórica incluía la condición de normalidad. Con el tiempo, la comunidad se dio cuenta de que esa condición no es estrictamente necesaria, la teoría es igualmente atractiva sin ella. Siguiendo esta tendencia, González y Teissier lograron extender el teorema \ref{GS forever} $(i)$ de González-Sprinberg quitando la hipótesis de normalidad. 
\end{rem}

Exhaustos y con la moral por los suelos, se viene un nuevo silencio de ocho años, interrumpido únicamente por una breve incursión al caso de curvas tóricas en la tesis de licenciatura de Daniel Green \cite{G,DG}.

\subsection{El regreso de la característica positiva (2020-2024)}

Y entonces, ocurrió un milagro. No al tercer día, sino cuarenta y cinco años después, la explosión de Nash resucitaría sobre campos de característica positiva. 

Retrocediendo en el tiempo, recordemos que en 1975 Nobile calculó un ejemplo que mostraba que la explosión de Nash no resuelve singularidades en característica positiva (ejemplo \ref{Nob ejem}). Coloquialmente hablando, ese ejemplo no es normal, ¡en el sentido más matemático posible! 

En una colaboración con Luis Núñez Betancourt logramos demostrar que añadiendo la hipótesis de normalidad a la variedad algebraica, el fenómeno de Nobile desaparecía. 

\begin{teo}\cite[Theorem 3.10]{DN1}
Sea $X$ una variedad algebraica normal sobre un campo de característica prima. La explosión de Nash $\nu:X^*\to X$ es un isomorfismo si y sólo si la variedad $X$ no tiene singularidades.
\end{teo}

Como es de esperarse, la explosión de Nash no preserva la normalidad por lo que debemos inmediatamente normalizar. ¡Y estamos de vuelta con la ENN! 

Inmediatamente nos pusimos a investigar las propiedades de la ENN en característica positiva y los resultados no se hicieron esperar. Siendo terreno inexplorado, decidimos empezar por caminos ya recorridos en característica cero. Al equipo Duarte-Núñez se une Jack Jeffries para lograr los primeros resultados. Ese trabajo generalizó a campos de característica positiva todo lo que se sabía en característica cero sobre la ENN de variedades tóricas.

\begin{teo}\cite[Theorems 1.9, 2.5]{DJNB}\label{sup tor car p}
Sea $p>0$ la característica del campo base $\K$.
\begin{itemize}
\item[(i)] La explosión de Nash de una variedad tórica coincide con la explosión del ideal logarítmico jacobiano módulo $p$.
\item[(ii)] La iteración de la ENN resuelve las singularidades de las superficies tóricas normales sobre $\K$.
\end{itemize}
\end{teo}

\begin{rem}\label{log jac 3}
\textit{Ideal logarítmico jacobiano (tercera parte).} Retomando las palabras de Spivakovsky, nuestro trabajo en característica $p$ generalizó, y fue inspirado por, el trabajo de González-Sprinberg, Lejeune-Jalabert, Reguera, González y Teissier. Con esto, las herramientas teóricas para el estudio combinatorio de la explosión de Nash y su versión normalizada de variedades tóricas, en cualquier característica, estaban listas.
\end{rem}

Con la búsqueda de nuevas familias para explorar la ENN en característica positiva, da inicio una nueva colaboración con Thaís Dalbelo y Maria Aparecida Soares Ruas que nos lleva a estudiar aquello que hicieron Ebelling y Gusein-Zade sobre variedades determinantales genéricas. 

\begin{teo}\cite[Theorem 4.3]{DDS}\label{det gen car p}
Sea $\K$ un campo de característica prima. La explosión de Nash de las 2-variedades determinantales genéricas sobre $\K$ es no singular.
\end{teo}

En este momento el entusiasmo estaba desbordado: ¡la explosión de Nash funciona en característica prima!

Ahora bien, podemos resumir en una línea la idea de la demostración tanto para las superficies tóricas como para las 2-variedades determinantales genéricas: la combinatoria involucrada en esos casos no depende de la característica del campo. También sabemos que eso no es cierto para otras variedades tóricas (ver \cite[Section 3]{DJNB}). Esto nos lleva a la pregunta obligada: 
\begin{center}
\textit{¿Bajo qué condiciones la combinatoria de la explosión de Nash de \\ una variedad tórica es independiente de la característica? }
\end{center}
Esta pregunta marcaría el inicio del fin.

\subsection{El contraejemplo (2024)}

Tratando de resolver la pregunta planteada previamente, Maximiliano Leyton, chileno de sabiduría serena, entra en escena. En un par de visitas al Centro de Ciencias Matemáticas, UNAM Campus Morelia, iniciamos una discusión al respecto. Múltiples ideas tuvieron lugar, el proyecto estaba en marcha. Tiempo después, Maximiliano me invitó a la Universidad de Talca, en Chile, a continuar con la colaboración. Las casualidades de la vida quisieron que en esa misma visita un encuentro de geometría algebraica estuviera sucediendo en Talca. 

El encuentro en cuestión es conocido como \href{https://sites.google.com/view/agrega0/home}{AGREGA} (Agrupando Equipos en Geometría Algebraica). La dinámica es la siguiente:  tres conferencistas presentan algún problema y, entre todo el público, se intentarán resolver en la semana que dura el evento. Yo fui uno de esos tres conferencistas en aquella ocasión. Tras largas y, debo decir, sádicas sesiones, dos personas más se unieron al equipo que habíamos iniciado Max y yo: Federico Castillo, mago colombiano de la programación, y Álvaro Liendo, chileno de inquebrantable determinación. El equipo estaba completo.

La discusión inició con la pregunta de la independencia de la característica. Las ideas que teníamos Max y yo se basaban en ejemplos sencillos en dimensión tres, calculables a mano. Pero era necesario iterar y los cálculos se complicaban. Era momento de llevar a terrenos computacionales la descripción combinatoria discutida en las observaciones \ref{log jac 1}, \ref{log jac 2} y \ref{log jac 3} sobre el ideal logarítmico jacobiano. Federico y Álvaro, expertos programadores, implementaron los algoritmos correspondientes en SageMath \cite{Sage}. Eso dio paso a una intensa exploración computacional, meses de cálculos, depuración de los algoritmos, búsqueda de invariantes que mejoraran en cada iteración. Todo en vano: cada ejemplo calculado concluía en una resolución de singularidades sin ninguna razón u orden aparente. Esto continuó por varios meses hasta que un brinco en la dimensión hizo la diferencia.

Por mucho tiempo los cálculos se concentraron en variedades tóricas de dimensión tres. Por la frustración que nos generaba no entender por qué el algoritmo parecía funcionar en ese caso, decidimos subir a dimensión cuatro. Este paso fue acompañado de una observación de Federico. Él, ajeno a las tribulaciones de los singularistas, sugirió estudiar una familia que suele ser fuente de contraejemplos en su área de especialidad, a saber, las variedades tóricas definidas por los conos Reeves. No pasó mucho tiempo cuando el programa entró en un bucle (ver figura 4). 
\begin{figure}[ht]\label{whats}
\begin{center}
\includegraphics[scale=0.3]{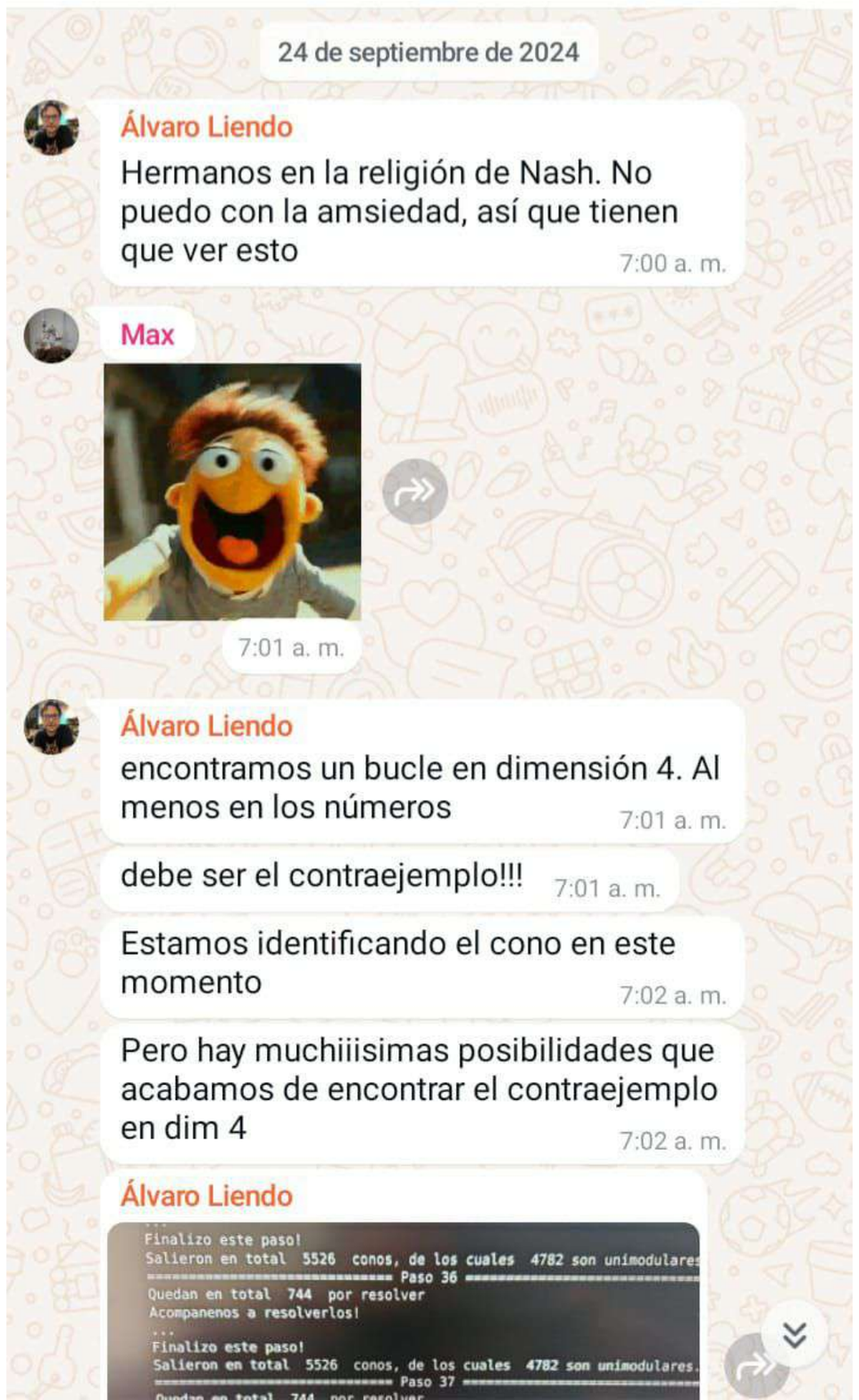}
\caption{El día en que se encontró el contraejemplo.}
\end{center}
\end{figure}

Las alarmas sonaron a todo vuelo, ¡el contraejemplo estaba a la vista! Unas horas más tarde los programas estaban dando a luz a la siguiente criatura:
\begin{align}
X=\V(&x_1x_6-x_3x_5, x_4x_6-x_2x_7, x_7^2-x_2x_4x_5,\notag\\
&x_6x_7-x_2^2x_5,x_3x_7-x_1x_2^2, x_3x_4x_5-x_1x_2x_7)\subset\C^7.\notag
\end{align}

La variedad tórica $X\subset\C^7$ tiene la siguiente particularidad: la explosión de Nash de $X$, que vive dentro de un espacio proyectivo, tiene una cubierta dada por variedades tóricas afines y, en una terrible burla del destino, una de estas variedades es isomorfa a $X$. En otras palabras, la explosión de Nash de $X$ contiene una copia idéntica a $X$; iterar el procedimiento no resolverá sus singularidades. Sucesiva exploración computacional nos llevó a encontrar contraejemplos también sobre campos de característica positiva. 

Estos ejemplos dieron lugar al artículo \cite{CDLL}, que fue aceptado para su publicación en \textit{Annals of Mathematics}.

\subsection{¿El fin? (2025-$\cdots$)}

En vista del contraejemplo que muestra que el método de Nash no resuelve singularidades en general, es obligado preguntarse: ¿la explosión de Nash ha llegado a su fin? ¡De ninguna manera! Paradójicamente, está gozando de una renovada vitalidad. 

En los últimos años han aparecido nuevos trabajos sobre las preguntas de Nash para variedades tóricas: dos tesis doctorales dedicadas a dimensión tres \cite{B,Su} y una tesis de maestría a dimensión dos \cite{C1}, dos artículos con respuestas afirmativas para nuevas familias en dimensión dos \cite{DS,C2}, nueva evidencia positiva en dimensión tres \cite{CDLL2}, un nuevo contraejemplo en dimensión tres \cite{CDLL3} y, finalmente, un artículo con la implementación computacional de los algoritmos \cite{CDLL4}. En otra dirección, un artículo resume las diferentes versiones del teorema de Nobile \cite{No2} y más recientemente se obtuvo una versión real y biLipschitz de ese mismo teorema \cite{Sam}.

\medskip

Paciente lectora,  atento lector: esta historia no ha terminado.

\vspace{.3cm}
\noindent{\footnotesize \textsc {D. Duarte\\ Centro de Ciencias Matem\'aticas, UNAM Campus Morelia.} \\
E-mail: adduarte@matmor.unam.mx}

\end{document}